\documentclass[11pt]{article}
\usepackage{amssymb}
\usepackage{amsbsy}
\usepackage{amsmath}

\def\BE{\begin{equation}}
\def\EE{\end{equation}}
\def\BQ{\begin{quote}}
\def\EQ{\end{quote}}
\def\si{\sigma}
\def\hh{{\bf h}}
\def\qed{$\bullet$}

\def\2nx2n{\mbox{\tiny{$2n\times 2n$}}}

\def\n1xn1{\mbox{\tiny{$(n-1)\times (n-1)$}}}
\def\nx1{\mbox{\tiny{$n \times 1$}}}

\def\0{\mbox{\tiny{$0$}}}
\def\1{\mbox{\tiny{$1$}}}
\def\2{\mbox{\tiny{$2$}}}
\def\3{\mbox{\tiny{$3$}}}
\def\4{\mbox{\tiny{$4$}}}
\def\-{\mbox{\tiny{$-$}}}
\def\+{\mbox{\tiny{$+$}}}

\def\H{\mathbb{H}}

\def\R{\mathbb{R}}

\def\inv{{\mbox{\tiny $-1$}}}
\def\mbf{\mathbf}

\begin{document}

\bibliographystyle{plain}
\title{
Quaternionic Commutations}
%\thanks{Submitted to
%{\sl Electronic Journal of Linear Algebra} on 18 July 1999.}}
%Editors will write the exact dates.

\author{
Nir Cohen\thanks{
Department of Applied Mathematics, IMECC, University of Campinas,
CP 6065, 13081-970  Campinas (SP) Brazil
(nir@ime.unicamp.br, deleo@ime.unicamp.br)}
\and
Stefano De Leo\footnotemark[2] \and Gisele C. Ducati\thanks{
Department of Mathematics, Federal University of Parana,
CP 19081, 81531-990  Curitiba (PR) Brazil
(ducati@mat.ufpr.br) }}

%If the same address is for the first and second person,
%\footnotemark[2] should be used.
%
%If the same address is for the second and third person,
%\footnotemark[3] should be used.

%authors and running title to go on top of the page
\markboth{N.\ Cohen, S.\ De Leo, G. C. Ducati}
{Quaternionic Commutations}

\maketitle

\begin{abstract} {\it
Given $n$ quaternions we investigate the extent
of non-commutativity of their multiple products,
commutators and exponential products.}
\end{abstract}

%\begin{keywords}
%Quaternions, Matrices, Determinants.
%\end{keywords}
%\begin{AMS}
%15A09, 15A33 \end{AMS}

%%%%%%%%%%%%%%%%%%%%%%%%%%%%%%%%%%%%%%%%%%%%%%%%%%%%%%%%%%%%%%%%%%%%%%%%%%%%%%%
%                                SECTION I
%%%%%%%%%%%%%%%%%%%%%%%%%%%%%%%%%%%%%%%%%%%%%%%%%%%%%%%%%%%%%%%%%%%%%%%%%%%%%%%

%\subsection{Introduction}
\section{Introduction}

The field of quaternions, $\H$, is not commutative. However,
it has some weak commutation properties, which we wish to
point out in this letter.

First, given $n$ quaternions $q1$ through $qN$ all the
multicommutators of the form
$$C(q1,\cdots,qN;\si):=[q_{\si(1)},[q_{\si(2)},\cdots
,[q_{\si(n-1)},q_{\si(n)}]\cdots]]~,$$ parametrized by
the various permutations $\si\in S_n,$ are pairwise equal
up to a $\pm$ sign.

Secondly, all the multiproducts of the form \BE\label{n-prod}
P(q1,\cdots,qN;\si):=q_{\si(1)}q_{\si(2)}\cdots q_{\si(n)}\EE for
which $\si$ is cyclic (consists of a single cycle) are mutually
similar. Thus, the $n!$ multiproducts defined by $q1$ through $qN$
occupy at most $(n-1)!$ similarity classes.

Finally, the Campbell-Baker-Hausdorf formula for the product
of two non-commuting exponentials, which in general
contains an infinite sum of non-commuting words in $p$
and $q$, reduces over the quaternions to a simple 4-term sum.

Perhaps the most basic ``commutation property'' of quaternions is the
commutability under the norm sign. Recall that the norm (or absolute
value) $|q|$ of $q=q_{\0}+ q_{\1} i + q_{\2} j + q_{\3} k$ is defined by
$|q|^2= q \bar q = q_{\0}^{\2}+ q_{\1}^{\2} + q_{\2}^{\2} + q_{\3}^{\2}~.$ It
satisfies the two commutation identities
\BE\label{norm}
|pq|=|qp|=|q| \, |p|~~~~~\mbox{and}~~~~~|1 - pq|=|1 - qp|~.
\EE

%%%%%%%%%%%%%%%%%%%%%%%%%%%%%%%%%%%%%%%%%%%%%%%%%%%%%%%%%%%%%%%%%%%%%%%%%%%%%%%
%                                SECTION II
%%%%%%%%%%%%%%%%%%%%%%%%%%%%%%%%%%%%%%%%%%%%%%%%%%%%%%%%%%%%%%%%%%%%%%%%%%%%%%%

\section{Basic product and commutation formulas}

Our notation for quaternions is the standard one, see \cite{codel}.
The quaternion $q \in \H$ is represented over the reals as \[
q=q_{\0}+ i \, q_{\1} + j \, q_{\2}+ k \, q_{\3}~,
~~~~~~~~~~~q_m \in \R~.
\]
We shall use the more concise vector notation $q= q_{\0}+\hh \cdot
{\mbf q}$, where $\hh=(i,j,k)$ and ${\mbf
q}=(q_{\0},q_{\1},q_{\2})\in\R^3.$ Using the inner and outer products in
$\R^3$, denoted here respectively by ${\mbf a}\cdot{\mbf b}$ and
${\mbf a}\times {\mbf b}$, we may write a concrete formula for the
product of two quaternions: \BE \label{2-prod}
ab=(a_{\0}b_{\0}-{\mbf a} \cdot {\mbf b})+\hh\cdot(a_{\0}{\mbf
b}+b_{\0}{\mbf a}+ {\mbf a}\times {\mbf b}). \EE As for the
commutator, it follows from (\ref{2-prod}) that \BE\label{2-com}
[a,b]=2\hh\cdot({\mbf a}\times {\mbf b}). \EE It also follows from
(\ref{2-com}) that a quaternion $q$ is a commutator if and only if
$q_{0}:=Re[q]=0.$

The commutator formula in (\ref{2-com}) easily generalizes to the
$n$-commutators $C(q1,\cdots,qN;\si)$ defined in the Introduction,
obtaining the formula \BE\label{com} C(q1,\cdots,qN;\si)=2^m
sgn(\si) \hh\cdot({\mbf q}{\mbf 1}\times\cdots\times {\mbf q}{\mbf
N}). \EE Thus, the various commutators defined by $q1$ through
$qN$ are all equal up to a $\pm$ sign.

\section{Similarity}

%%%%%%%%%%%%%%%%%%%%%%%%%%%%%%%%%%%%%%%%%%%%%%%%%%%%%%%%%%%%%%%%%%%%%%%%%%%%%%%
%                                SECTION III
%%%%%%%%%%%%%%%%%%%%%%%%%%%%%%%%%%%%%%%%%%%%%%%%%%%%%%%%%%%%%%%%%%%%%%%%%%%%%%%

%\section*{Similarity of products}

Two quaternions $p$ and $q$ are called similar if
\[
q=s^\inv p \, s~,~~~~~~~ \mbox{for some}~s\in\H~.
\]
It can be seen that similarity of $p$ and $q$ amounts to the two
conditions $|p|=|q|$ and $Re[p]=Re[q].$ Now we ask whether the two
quaternion products $p:=P(q1,\cdots,qN;\si_1)$ and
$q:=P(q1,\cdots,qN;\si_2)$ are similar. The first condition
$|p|=|q|$ is guaranteed in view of (\ref{norm}). Thus, similarity
of $p$ and $q$ in (\ref{n-prod}) is reduced to the second
condition, $Re[q]=Re[p].$

The condition for
similarity of $2,3$- and $4$-products is given below.
The general case is still open.

\begin{quote}{\bf Lemma. }{\it
In the non-commutative field of quaternions:

(i) For all $p,q$ the products $pq$ and $qp$ are always similar;

(ii) If in (\ref{n-prod}) $q$ is obtained from $p$
by a primitive permutation $(\si(i)=i+k~mod(n))$
then $p$ and $q$ are similar;

(iii) $p=abc$ and $q=acb$ are similar if and only if
the vectors ${\mbf a},{\mbf b},{\mbf c} \in \R^3$ are linearly dependent.

(iv) $p = abcd$ and $q = adcb$ are similar if and only if
$~a_0 \alpha - b_0 \beta + c_0 \gamma - d_0 \delta = 0$,
 $\alpha, \beta,\gamma,\delta \in \R$.}

\end{quote}

\noindent{\bf Proof. }
(i) Assume that $p=ab$ and $q=ba$ for some $a,b\in\H$.
If $b=0$ then $p=q=0$ and there is nothing to prove.
Otherwise, the identity $bpb^\inv=q$ holds, implying that
$p$ and $q$ are similar.

(ii) Follows directly from (i).

(iii) By (ii), the six terms $P(a,b,c;\si)$ $(\si\in S_3)$
can occupy at most two similarity classes, represented by
$abc$ and $acb.$ The only requirement for similarity of these
two terms is that $Re(abc) = Re(acb).$
Direct calculation based on (\ref{2-prod}) gives
\begin{eqnarray}
Re(abc)=(ab)_{\0}c_{\0}-{\mbf (ab)} \cdot {\mbf c}= \hspace*{2cm}\nonumber\\
=a_{\0}b_{\0}c_{\0}-[a_{\0}({\mbf b}\cdot{\mbf c})+ b_{\0}({\mbf
a} \cdot {\mbf c})+ c_{\0}({\mbf a}.{\mbf b})]-b{\mbf a}\times
{\mbf b} \cdot {\mbf c}.
\end{eqnarray}
It follows that
$$Re[abc]-Re[acb]=2 \, ({\mbf a}\times {\mbf b}\cdot{\mbf c})=
2 \, \det[{\mbf{a,b,c}}],$$ proving the claim.
It is observed that for a generic triple $(a,b,c)$ the terms $abc$
and $acb$ are not similar.

(iv) The case $n=4$ has six different similarity classes,
represented by $abcd,$ $abdc,$ $acbd,$ $acdb,$ $adbc,$ $adcb.$
All six expressions have the same norm, hence one only need
compare their real parts.

Some cases can be reduced to $n=3$. For example, considering $ab$ a fixed
vector we have
\begin{eqnarray*}
Re[abcd]-Re[abdc]= 2 \,[ a_{\0} {\mbf b} + b_{\0} {\mbf a} + {\mbf a} \times
{\mbf b}] \times {\mbf d} \cdot {\mbf c}
\end{eqnarray*}
The vectorial part of $ab$ is
$a_{\0} {\mbf b} + b_{\0} {\mbf a} + {\mbf a} \times {\mbf b}$ and we have
the same result obtained for $n = 3$ as expected.
The same holds to $abcd$ and $adbc$.
In this case $bc$ is the fixed vector.

The pair $abcd$ and $acbd$ must be analysed. We find
\begin{eqnarray}
\label{case4}
Re[abcd]-Re[adcb]= -2 \left[ (a_{\0} {\mbf b} + b_{\0} {\mbf a})\cdot (
{\mbf c} \times {\mbf d}) + ({\mbf a} \times {\mbf b})\cdot (c_{\0} {\mbf d} +
d_{\0} {\mbf c}) \right]
\end{eqnarray}

Two cases should be analysed. If $span\{ {\mbf{a, b, c, d}} \} =
\R^2$, (\ref{case4}) is null, as is easy to check.
If $span\{ {\mbf{a, b, c, d}} \} = \mathbb{R}^3$ consider
\begin{equation}
\label{d}
\alpha {\mbf a} + \beta {\mbf b} + \gamma {\mbf c} + \delta {\mbf d} =
0~,~~~\alpha, \beta, \gamma, \delta \in \mathbb{R}~.
\end{equation}
Substituting (\ref{d}) into (\ref{case4}) yields
\[
a_0 \alpha - b_0 \beta + c_0 \gamma - d_0 \delta = 0~,
\]
which completes the proof. \qed

Remarks:

1. With item (ii) we may reduce the number of similarity classes
involving permutations of a given $n$-product from $n!$ to
$(n-1)!$ at most. Can this number be reduced further?

2. We did not investigate the extent to which our observations
generalize to other rings of interest. Note that (i,ii) will
hold in every integral domain.

3. Investigate the equality of norms for the various expressions
of the form $P(q1,\cdots,qN;\si)-t,$ $t\in\R.$ It follows from
(\ref{2-prod}) that cyclic permutation has no effect on the norm.
What is the geometric interpretation? What about non-cyclic
permutations?

4. We should search the conditions for true equality between two
multiproducts.

%5. Are $P(q1,\cdots,qN;\si)$ always similar when $span\{ {\mbf{q1,
%\cdots, qN}} \} = \R^2$?

\vspace*{.5cm}

%%%%%%%%%%%%%%%%%%%%%%%%%%%%%%%%%%%%%%%%%%%%%%%%%%%%%%%%%%%%%%%%%%%%%%%%%%%%%%%
%                                SECTION IV
%%%%%%%%%%%%%%%%%%%%%%%%%%%%%%%%%%%%%%%%%%%%%%%%%%%%%%%%%%%%%%%%%%%%%%%%%%%%%%%

\section{The Campbell-Baker-Hausdorff formula}

As is well known, the formula $Exp[p]Exp[q]=Exp[p+q]$ does not
hold in general, when $p$ and $q$ do not commute. The well-known
CBH (Campbell-Baker-Hausdorff) formula contains an additive
correction term within the right hand side exponential which
allows the equality to be restored. In the general case, the
correction term is a certain infinite sum involving non-commuting
words of arbitrary length made of the letters $p$ and $q$;
however, in the case of quaternions  this formula should be
closed.

A detailed discussion of the quaternionic CBH formula and its
closed form will appear in a forthcoming paper [Ref.: Quaternionic
exponentials, in preparation].

%%%%%%%%%%%%%%%%%%%%%%%%%%%%%%%%%%%%%%%%%%%%%%%%%%%%%%%%%%%%%%%%%%%%%%%%%%%%%%%
%                                SECTION V
%%%%%%%%%%%%%%%%%%%%%%%%%%%%%%%%%%%%%%%%%%%%%%%%%%%%%%%%%%%%%%%%%%%%%%%%%%%%%%%

\section{Exponential derivative}

Due the non commutativity of quaternions
\begin{equation}
\label{od} \{ \exp[{\psi(x)}] \}' = \psi'(x) \exp{[\psi(x)]}~.
\end{equation}
{\it does not} hold. Let us see how the expression (\ref{od}) is
changed. From the different ways to write a quaternion function we
have chosen the more convenient one for this section. In order to
find the derivative of the exponential function we shall consider
\[
\psi(x) = f(x) + I(x) g(x)
\]
where
\[
f(x) = \psi_{\0}(x)~,~~~g(x) = \sqrt{\psi_{\1}^{\2}(x) + \psi_{\2}^{\2}(x)
+ \psi_{\3}^{\2}(x)}
\]
are real functions of $x$ and
\[
I(x) = \frac{i \psi_{\1}(x) + j \psi_{\2}(x) + k \psi_{\3}(x)}
{\sqrt{\psi_{\1}^{\2}(x) + \psi_{\2}^{\2}(x) + \psi_{\3}^{\2}(x)}}
\]
with $|I(x)|^{\2} = -1$. Note that
\begin{equation}
\label{anti}
\{ I(x), I'(x) \} := I(x) I'(x) + I'(x) I(x) = 0
\end{equation}
where $'$ denotes the first derivative.
Consider ${\mathcal F'}(x) = f'(x) + I(x) g'(x)$.
Omitting the $x$ variable we have
\[
\psi = f + I g~~~\to~~~\psi' = f' + I g' + I'g =
{\mathcal F'} + I' g~.
\]
%Now we shall show that
%\[
%e^{\psi} = \exp[f + I g] = e^{f} \left( \cos g +
%I \sin g  \right)
%\]
Using power series expansion of the exponential function~\cite{HAM} we find
\BE
\label{exp}
e^{\psi}  =  e^{f} e^{I g} = e^{f} \sum_{n=0}^{\infty}
\frac{[Ig]^n}{n!} = e^f \{ \cos g + I \sin g \}
%  & = & e^f \left\{ \sum_{k=0}^{\infty} \frac{I^{2k} g^{2k}}{(2k)!} +
%\sum_{k=0}^{\infty} \frac{I^{2k+1} g^{2k+1}}{(2k+1)!} \right\} \nonumber\\
% & = & e^f \left\{ \sum_{k=0}^{\infty} \frac{(-1)^{k} g^{2k}}{(2k)!} +
%I \, \sum_{k=0}^{\infty} \frac{(-1)^{k} g^{2k+1}}{(2k+1)!} \right\}
%\nonumber\\
% & = & e^f \{ \cos g + I \sin g \}
\EE
Thus, the exponential function derivative, using (\ref{exp}) and the property
(\ref{anti}), is given by
\begin{eqnarray*}
[ e^{\psi} ]' & = & {\mathcal F'} e^\psi + I' e^f
                       \sin g \nonumber\\
                 & = & [ {\mathcal F'} + I'g - I'g ] e^\psi
                      + I' e^f \sin g  \nonumber\\
                 & = & \psi' e^\psi - I' \left\{
                       g e^\psi - e^f \sin g  \right\}
\end{eqnarray*}

\end{document}